\documentclass[12pt]{article}

\usepackage{amssymb}

\newtheorem{theorem}{Theorem}[section]
\newtheorem{corollary}[theorem]{Corollary}

\newtheorem{proposition}[theorem]{Proposition}

\begin{document}
	
\title{On solution manifolds for algebraic-delay systems} 
	
\author{ Hans-Otto Walther\\Mathematisches Institut\\Universit\"at Gie{\ss}en}

\date{\small{{\tt Hans-Otto.Walther@math.uni-giessen.de}}}

\maketitle

\begin{abstract}
\noindent
Differential equations with state-dependent delays define a semiflow of continuously differentiable solution operators in general only on an associated submanifold of the Banach space $C^1([-h,0],\mathbb{R}^n)$. We extend a recent result on simplicity of these {\it solution manifolds} to systems where the delay is given by the state only implicitly in an extra equation. Such algebraic-delay systems arise from various applications. 
\end{abstract}
	
\bigskip
	
	
\medskip


\bigskip
 
\section{Introduction}

For $h>0$ and $n\in\mathbb{N}$ let $C_n$ and $C^1_n$ denote the Banach spaces of continuous and of continuously differentiable maps $[-h,0]\to\mathbb{R}^n$, respectively. In case $n=1$ we abbreviate $C=C_1$, $C^1=C^1_1$. For a map $x:[-h,t_e)\to\mathbb{R}^n$ and $0\le t<t_e\le\infty$ the segment $x_t:[-h,0]\to\mathbb{R}^n$ is defined by $x_t(s)=x(t+s)$.

Delay differential equations in the general form
\begin{equation}
x'(t)=f(x_t) 
\end{equation}
with a continuously differentiable map 
$f:C^1_n\supset U\to\mathbb{R}^n$ which satisfies an additional smoothness condition define a continuous semiflow of continuously differentiable solution operators on the {\it solution manifold}
$$
X_f=\{\phi\in U:\phi'(0)=f(\phi)\}
$$
which in case $X_f\neq\emptyset$ is continuously differentiable and has codimension $n$ in $C^1_n$ \cite{W1,HKWW}. The additional smoothness property just mentioned is that the derivatives $Df(\phi):C^1_n\to\mathbb{R}^n$, $\phi\in U$, have linear extensions $D_ef(\phi)$ to the space $C_n$ so that the map
$$
U\times C_n\ni(\phi,\chi)\mapsto D_ef(\phi)\chi\in\mathbb{R}^n
$$
is continuous. See \cite{M-PNP} for the first version of this extension property, which holds true for a large class of concrete equations with state-dependent delays when written in the general form of Eq. (1).
Let us recall that for differential equations with state-dependent delay the initial value problem for data in $C_n$ or $C^1_n$ is in general not well-posed \cite{Wi,HKWW}. 

We are interested in the nature of solution manifolds, whether they are simple as a graph with respect to a direct sum decomposition of the space $C^1_n$, or whether they can be more complicated. From the defining equation $\phi'(0)=f(\phi)$ it is only obvious that for $f=0$ the trivial solution manifold 
$$
X_0=\{\chi\in C^1_n:\chi'(0)=0\}
$$
is a closed subspace of codimension $n$ in $C^1_n$. Results in \cite{KR,W6} guarantee that for $D_ef$ bounded or for delays bounded away from zero (formulated as a condition on $f$) solution manifolds can be written as graphs over $X_0$ while a counterexample in \cite{W6} shows that in general this is impossible; in general solution manifolds do not admit a  graph representation with respect to any direct sum decomposition of the space $C^1_n$. 

 In \cite{KW2} we proved for a  class of systems 
\begin{equation}
x'(t)=g(x(t-d_1(Lx_t)),\ldots,x(t-d_{k}(Lx_t))) 
\end{equation}
with discrete state-dependent delays that the associated solution manifolds can be transformed into open subsets of the space $X_0$, by diffeomorphisms which leave  points of $X_f\cap X_0$ fixed. This means that they are {\it almost graphs} over $X_0$,  in the terminology of \cite{W4,W6,KW2}. An example of an almost graph in the plane is the unit circle without its uppermost point.

The assumptions on Eq. (2) in \cite{KW2} are that $g:\mathbb{R}^{kn}\supset V\to\mathbb{R}^n$ is continuously differentiable, that $L:C_n\to F$ is a continuous linear map into a finite-dimensional normed real vectorspace $F$, that $W\subset F$ is open, that the delay functions $d_{\kappa}:W\to[0,h]\subset\mathbb{R}$, $\kappa=1,\ldots,k$, are continuously differentiable, and, of course, that there exist $\phi\in C^1_n$ with $L\phi\in W$ so that the vector with the components $v_{\iota}=\phi_j(-d_{\kappa}(L\phi))$, for $j\in\{1,\ldots,n\}$ and $\kappa\in\{1,\ldots,k\}$ determined by $\iota=(\kappa-1)n+j$, belongs to $V$.

As a motivation for the form $d_{\kappa}(L\phi)$ of the delays in Eq. (2) one may think of $L\phi$ approximating $\phi\in C_n$ in a finite-dimensional subspace $F\subset C_n$.

In several
equations arising from applications delays are not always explicit as in Eq. (2) but are given implicitly by an additional equation involving the delay and the state of the system - see for example a recent model for protein synthesis \cite{GHMWW}, or models for position control by echo \cite{W2,W3}. A general form of such algebraic-delay systems is
 \begin{eqnarray}
x'(t)=G(r(t),x_t),\\ 
0 = \Delta(r(t),x_t) 
\end{eqnarray}
Some basic theory for Eqs. (3,4) was developped in \cite{W4}, for 
continuously differentiable maps $G:\mathbb{R}^k\times C^1_n\supset {\mathcal U}\to\mathbb{R}^n$ and $\Delta:{\mathcal U}\to\mathbb{R}^k$ which have the following additional smoothness property:

(He)  {\it All derivatives $DG(s,\phi):\mathbb{R}^k\times C^1_n\to\mathbb{R}^n$ and $D_2\Delta(s,\phi):C^1_n\to\mathbb{R}^k,$ $(s,\phi)\in {\mathcal U}$,
have linear extensions $D_eG(s,\phi)$ and $D_{2,e}\Delta(s,\phi)$ to $\mathbb{R}^k\times C_n$ and to $C_n$, respectively, so that both maps
$$
{\mathcal U}\times\mathbb{R}^k\times C_n\ni(s,\phi,p,\chi)\mapsto
D_eG(s,\phi)(p,\chi)\in\mathbb{R}^n
$$
and
$$
{\mathcal U}
\times C_n\ni(s,\phi,\chi)\mapsto D_{2,e}\Delta(s,\phi)\chi\in\mathbb{R}^k
$$
are continuous.}

Moreover it is assumed in \cite{W4} that $r\in(-h,0)^k$  for each $(r,\phi)\in {\mathcal U}$. The main result  of \cite{W4} says that the set
$$
M_{G,\Delta}=\{(s,\phi)\in
U:\phi'(0)=G(s,\phi),0=\Delta(s,\phi),\det\,D_1\Delta(s,\phi)\neq0\},
$$
if non-empty, is a continuously differentiable
submanifold of codimension $k+n$ in $\mathbb{R}^k\times C^1_n$ on which
Eqs. (3,4) define a continuous semiflow of continuously differentiable solution operators. The assumption concerning $r\in (-h,0)^k$ restricts applications to systems whose delays have no zeros. In order to remove this restriction one may replace the previous condition by the weaker hypothesis that there is an open interval $I\supset[-h,0]$ so that
$r\in I^k$ for all $(r,\phi)\in{\mathcal U}$. Then the main result of \cite{W4} remains valid, with the same proof, in which only 4 formulations must be adjusted. Section 6 below contains the details.

We call $M_{G,\Delta}$ the solution manifold associated with Eqs. (3,4). 

Let us recall here the notion of a solution of  Eqs. (3.4) on an interval $[-h,t_e)$, $0<t_e\le\infty$, which is a pair $(x,r)$ of a continuously differentiable map $x:[-h,t_e)\to\mathbb{R}^n$ together with a continuous map $r:[0,t_e)\to\mathbb{R}^k$ for which Eqs. (3,4) hold for $0\le t<t_e$.

Incidentally let us recall the following example from \cite{W31,W4} : The differential equation for a feedback system which reacts to its present state $x(t)\in\mathbb{R}$ only after a delay $d(x(t))\ge0$ reads
\begin{equation}
x'(t+d(x(t)))=f(x(t)), 
\end{equation}
with functions $d:\mathbb{R}\to[0,\infty)$ and $f:\mathbb{R}\to\mathbb{R}$. The attempt to rewrite Eq. (5) as a delay differential equation in the familiar form $x'(s)=\ldots$ yields
\begin{equation}
x'(s)=f(x(s+r(s))) 
\end{equation}
with $s=t+d(x(t))$ and $r(s)=t-s\le0$, which implies
$$
r(s)=t-s=-d(x(t))=-d(x(s+r(s))),
$$
hence
\begin{equation}
0=d(x(s+r(s)))+r(s)=d(x_s(r(s)))+r(s). 
\end{equation}
The algebraic-delay system (6,7) has the form (3,4) with $G(r,\phi)=f(\phi(r))$ and $\Delta(r,\phi)=d(\phi(r))+r$.

In the sequel we study the solution manifold $M_{G,\Delta}$ for
algebraic-delay systems (3,4) which generalize the system (6,7) as well as Eq. (2). 

The description of these systems begins with the choice of an open interval $I\supset[-h,0]$ and of a larger compact interval $J\subset[-2h,h]$, and with the {\it odd extension map} $E:C\to C(J,\mathbb{R})$ defined by
\begin{eqnarray}
E\phi(t) & = & \phi(t)\quad\mbox{on}\quad[-h,0],\nonumber\\ 
E\phi(t) & = & 2\phi(-h)-\phi(-t-2h)\quad\mbox{for}\quad J\ni t\le-h,\nonumber\\ 
E\phi(t) & = & 2\phi(0)-\phi(-t)\quad\mbox{for}\quad0\le t\in J.\nonumber
\end{eqnarray}
The map $E$ is linear and continuous and induces a continuous linear map $C^1\stackrel{E}{\to}C^1(J,\mathbb{R})$. Both maps have norm 3. Let $k\in\mathbb{N}$ be given. It is convenient to associate with $(r,\phi)\in J^k\times C_n$ the vector $\widehat{(r,\phi)}=v\in\mathbb{R}^{kn}$ whose components are
$$
v_{\iota}=E\phi_j(r_{\kappa})
$$
with  $\kappa\in\{1,\ldots,k\}$ and $j\in\{1,\ldots,n\}$ determined by $\iota=(\kappa-1)n+j\in\{1,\ldots,kn\}$.

Let $g:\mathbb{R}^{kn}\supset V\to\mathbb{R}^n$ and $W\subset F$ be given as in Eq. (2). We consider Eqs. (3,4) for
\begin{eqnarray}
G(r,\phi) & = & g(\widehat{(r,\phi)}),\\ 
\Delta(r,\phi) & = & \delta(r,Q(r,\phi)) 
\end{eqnarray}
with continuous maps $\delta:J^k\times W\to\mathbb{R}^k$ and
$Q:J^k\times C_n\to F$ so that each map $Q(r,\cdot):C_n\to F$, $r\in J$, is linear. The domain of $G$ and $\Delta$ in (8,9) is
$$
{\mathcal U}=\{(r,\phi)\in I^k\times C^1_n:\widehat{(r,\phi)}\in V,Q(r,\phi)\in W\};
$$
in order that ${\mathcal U}$ be nonempty we assume that for some  $(r,\phi)\in I^k\times C^1_n$, $\widehat{(r,\phi)}\in V$ and $Q(r,\phi)\in W$. Concerning smoothness we assume that the restrictions of $\delta$ to $I^k\times W$ and of $Q$ to $I^k\times C^1_n$ are continuously differentiable. A final hypothesis concerns the maps 
$$
Q_j:J^k\times C\ni(r,\phi)\mapsto Q(r,\phi\cdot e_j)\in F,\quad j=1,\ldots,n.
$$
Here the product $\phi\cdot e_j\in C_n$, with the unit vector $e_j\in\mathbb{R}^n$, $e_{j,\nu}=\delta_{j\nu}$,  is defined by
its components $(\phi\cdot e_j)_j=\phi\in C$ and
$(\phi\cdot e_j)_{\nu}=0\in C$ for $\nu\neq j$ in
$\{1,\ldots,n\}$. 

We require that for every $j\in\{1,\ldots,n\}$ the map $q=Q_j$ satisfies
the following:

\noindent
(Hq) {\it The range $F_q\subset F$ of $q(r,\cdot)$ is constant (independent of $r\in J^k$), and either $F_q=0$ or there are $\beta_1,\ldots,\beta_d$, $d=\dim\,F_q$, in $C$ so that for every $r\in J^k$ the vectors $q(r,\beta_m)$, $m=1,\ldots,d$, form a basis of $F_q$.
}\\

If in condition (Hq) we have 
$F_q\neq0$ then it follows easily that $\beta_1,\ldots,\beta_d$ are linearly independent and that the space $\sum_{m=1}^d\mathbb{R}\beta_m$ is complementary in $C$ for each nullspace $q(r,\cdot)^{-1}(0)$, $r\in J^k$.

As condition (Hq) looks rather restrictive we give an example:
Select components of $\phi\in C_n$ by an injective map $\nu:\{1,\ldots,d\}\to\{1,\ldots,n\}$ and components of $r\in J^k$ by a map $\kappa:\{1,\ldots,d\}\to\{1,\ldots,k\}$ and define $Q:J^k\times C_n\to\mathbb{R}^d$ by
$$
Q(r,\phi)=(E\phi_{\nu(1)}(r_{\kappa(1)}),\ldots,E\phi_{\nu(d)}(r_{\kappa(d)})).
$$ 
Then each map $Q(r,\cdot)$ is linear. Using continuity of the evaluation map $ev:C(J,\mathbb{R})\times J\ni(s,\phi)\mapsto\phi(s)\in\mathbb{R}$  and continuous differentiability of the map
$ev_1:C^1(J,\mathbb{R})\times I\stackrel{ev}{\to}\mathbb{R}$ one sees that $Q$ is smooth as required. In order to verify property (Hq)
let $j\in\{1,\ldots,n\}$ be given. For $r\in J^k$, $\phi\in C$, and $m\in\{1,\ldots,d\}$, the $m-$th component of $w=q(r,\phi)=Q_j(r,\phi)=Q(r,\phi\cdot e_j)$ is  $E((\phi\cdot e_j)_{\nu(m)})(r_{\kappa(m)})$. In case $\nu(m)\neq j$ we get $w_m=0$ while in case $\nu(m)=j$, $w_m=E\phi(r_{\kappa(m)})$. Hence $Q_j(r,\phi)=0$ in case $j\notin\nu(\{1,\ldots,d\})$ while in case $j\in\nu(\{1,\ldots,d\})$ there is a unique $m\in\{1,\ldots,d\}$ with $\nu(m)=j$, the $m$-th component of $Q_j(r,\phi)$ is $E\phi(r_{\kappa(m)})$, and all other components are zero. It follows that either $Q_j(r,\cdot)=0$ for all $r\in J^k$, or there exists $m \in\{1,\ldots,d\}$ such that for every $r\in J^k$ we have $Q_j(\{r\}\times C)=\mathbb{R}\,e_m$, with the unit vector $e_m$ from the canonical basis of $\mathbb{R}^d$. In the last case, consider the constant function $\mathbf{1}:[-h,0]\ni t\mapsto1\in\mathbb{R}$ and notice that for every $r\in J^k$, $q(r,\mathbf{1})=Q_j(r,\mathbf{1})=e_m$ 
spans the range $F_q$ of $q(r,\cdot)=Q_j(r,\cdot)$.

Proposition 1.1 below guarantees that the system (3,4) with $G,\Delta$ specified by Eqs. (8,9) satisfies the smoothness hypotheses which in case $M_{G,\Delta}\neq\emptyset$ yield a continuous semiflow of continuously differentiable solution operators on the submanifold $M_{G,\Delta}$ of codimension $k+n$ in $\mathbb{R}^k\times C^1_n$.

The main result of the present paper is stated in Proposition 4.3 and Corollary 4.4 below. It guarantees that for the system  (3,4) with $G,\Delta$ specified by Eqs. (8,9) there is  a diffeomorphism $T:{\mathcal U}\to\mathbb{R}^k\times C^1_n$ onto an open subset ${\mathcal O}\subset I^k\times C^1_n\subset\mathbb{R}^k\times C^1_n$ which takes $M_{G,\Delta}$ to the subset of ${\mathcal O}$ given by 
$$
\psi'(0)=0,\quad\Delta(r,\psi)=0,\quad\det\,D_2\Delta(T^{-1}(r,\psi))\neq0.
$$
In case $M_{G,\Delta}\neq\emptyset$ the image $T(M_{G,\Delta})$ is a continuously differentiable submanifold of codimension $k+n$ in $\mathbb{R}^k\times C^1_n$ which is contained in the open strip $I^k\times X_0$ in the linear subspace $\mathbb{R}^k\times X_0$ of codimension $n$, and which is defined, apart from a regularity condition, by only the algebraic equation (4).

Our result also sheds more light on the solution manifold $X_f$ of Eq. (2), which is equivalent to a system (3,4,8,9), with 
$Q(r,\phi)=L\phi$ independent of $r$ and with $\delta(r,w)=d(w)+r$, $d=(d_1,\ldots,d_k)$. For this system we know $M_{G,\Delta}\neq\emptyset$, according to Proposition 1.2 below. Corollary 5.1 says that the solution manifold $M_{G,\Delta}$ is diffeomorphic to a graph in $I^k\times X_0$ given by $r=-d(L\phi)$  - whereas according to \cite{KW2} the solution manifold $X_f\subset C^1_n$ associated with Eq. (2) is an almost graph over $X_0$ and has in general no graph representation with respect to any direct sum decomposition of $C^1_n$.

The proof of Proposition 4.3 and Corollary 4.4 is by a modification of
the approach in \cite{KW2}. The diffeomorphism $T$ has the form
$T(r,\phi)=(r,A(r,\phi))$ with a parameterized version $A:\mathbb{R}^k\times C^1_n\supset {\mathcal U}\to C^1_n$ of a map from \cite{KW2} which at each point $\phi\in X_f$  is given by a projection onto the trivial solution manifold $X_0$ along a space which is complementary to both $X_0$ and to the tangent space $T_{\phi}X_f$. In the present paper we also need spaces which are complementary to $X_0$ in $C^1_n$. The construction of these spaces is prepared in Section 2 below. 

{\bf Notation, preliminaries.} On $\mathbb{R}^n$ we use the Euclidean norm. The vectors $c\in\mathbb{R}^n$ are columns, occasionally written as $(c_1,\ldots,c_n)^{tr}$, or even without the upper index when they appear as arguments of maps. 
The relation $A\subset\subset B$ for subsets of $\mathbb{R}^n$ means that the closure of $A$ is compact and contained in $B$.
Derivatives and partial derivatives indicated by capitals $D$ and $D_{\sigma}$, respectively, are continuous linear maps. For differentiable maps $\phi$ with domain in $\mathbb{R}$, $\phi'(t)=D\phi(t)1$, and for the  components of a differentiable map $m$ from an open subset of $\mathbb{R}^{s}$ into $\mathbb{R}^p$,   $\partial_{\sigma}m_{\pi}(x)=D_{\sigma}m_{\pi}(x)1$. 

The norm on the vectorspace $C(K,\mathbb{R}^n)$ of continuous maps from a compact interval $K$ into $\mathbb{R}^n$ is given by $|\phi|=\max_{t\in K}|\phi(t)|$. On the space $C^1(K,\mathbb{R}^n)$ of continuously differentiable maps $K\to\mathbb{R}^n$ the norm is given by\\ 
$|\phi|=\max_{t\in K}|\phi(t)|+\max_{t\in K}|\phi'(t)|$.

The odd extension map $E:C\to C(J,\mathbb{R})$ preserves constant functions in the sense that for every $c\in\mathbb{R},\quad E(c\mathbb{1})(t)=c$ on all of $J$.

The  differentiation map $\partial:C^1(K,\mathbb{R})\ni\phi\mapsto\phi'\in C(K,\mathbb{R})$ is linear and continuous.

For the map $ev_1:C^1(J,\mathbb{R})\times I\stackrel{ev}{\to}\mathbb{R}$ we have $D\,ev_1(\phi,t)(\chi,s)=\chi(t)+s\,\phi'(t)$.

\begin{proposition}
The domain ${\mathcal U}$ of $G$ and $\Delta$ given by (8,9) is open, both maps are continuously differentiable, and they have the smoothness property (He).
\end{proposition}

{\bf Proof.} 1.  On ${\mathcal U}$. The maps $J^k\times C_n\ni(r,\phi)\mapsto\widehat{(r,\phi)}\in\mathbb{R}^{kn}$ and $J^k\times C_n\ni(r,\phi)\mapsto Q(r,\phi)\in F$ are continuous. Using the inclusion maps $C^1_n\hookrightarrow C_n$ and $I\hookrightarrow J$ one sees that also the restrictions of the previous maps to $I^k\times C^1_n$ are continuous. So the preimages of the open sets $V$ and $W$ under the latter maps are open, hence ${\mathcal U}$ as the intersection of the open preimages is open.

2. On $G$.
Consider a component $G_{\nu}$, $\nu\in\{1,\ldots,n\}$, of $G$.
For every $(r,\phi)\in {\mathcal U}$,
$$
G_{\nu}(r,\phi)=g_{\nu}(h_{11}(r,\phi),\ldots,h_{n1}(r,\phi);\ldots;h_{1k}(r,\phi),\ldots,h_{nk}(r,\phi))
$$
with $h_{j\kappa}:{\mathcal U}\to\mathbb{R}$ given by $h_{j\kappa}(r,\phi)=E\phi_j(r_{\kappa})=ev_1(E\phi_j,r_{\kappa})$.
Using that coordinate projections are linear and continuous we obtain from the chain rule that each map $h_{j\kappa}$ is continuously differentiable with
$$
Dh_{j\kappa}(r,\phi)(s,\psi)=(E\psi_j)(r_{\kappa})+s_{\kappa}(E\phi_j)'(r_{\kappa})=ev(E\psi_j,r_{\kappa})+s_{\kappa}ev(\partial\,E\phi_j,r_{\kappa})
$$
for each $(r,\phi)\in {\mathcal U}\subset I^k\times C^1_n$ and all $(s,\psi)\in\mathbb{R}^k\times C^1_n$. It follows that $G_{\nu}$ is continuously differentiable with
\begin{eqnarray}
DG_{\nu}(r,\phi)(s,\psi) & = & \sum_{\kappa=1}^k\sum_{j=1}^n\partial_{(\kappa-1)n+j}g_{\nu}(\widehat{(r,\phi)})[Dh_{j\kappa}(r,\phi)(s,\psi)]\nonumber\\
& = &  \sum_{\kappa=1}^k\sum_{j=1}^n\partial_{(\kappa-1)n+j}g_{\nu}(\widehat{(r,\phi)})[ev(E\psi_j,r_{\kappa})+s_{\kappa}ev(\partial\,E\phi_j,r_{\kappa})]\nonumber
\end{eqnarray}
for each $(r,\phi)\in {\mathcal U}\subset I^k\times C^1_n$ and all $(s,\psi)\in\mathbb{R}^k\times C^1_n$. We infer that $G$ is continuously differentiable. 

The last expression for $DG_{\nu}(r,\phi)(s,\psi)$ also defines  linear maps $D_eG_{\nu}(r,\phi):\mathbb{R}^k\times C_n\to\mathbb{R}$, for every $\nu\in\{1,\ldots,n\}$ and for each $(r,\phi)\in {\mathcal U}$. Using  continuity of the differentiation map $\partial:C^1(J,\mathbb{R})\to C(J,\mathbb{R})$ and continuity of the evaluation map $ev:C(J,\mathbb{R})\times J\to\mathbb{R}$ we observe that every map
$$
{\mathcal U}\times\mathbb{R}^k\times C_n\ni(r,\phi,s,\chi)\mapsto D_eG_{\nu}(r,\phi)(s,\chi)\in\mathbb{R},\quad\nu\in\{1,\ldots,n\},
$$
is continuous. For every $(r,\phi)\in {\mathcal U}$ we define a linear map $D_eG(r,\phi):\mathbb{R}^k\times C_n\to\mathbb{R}^n$  by
$$
[D_eG(r,\phi)(s,\chi)]_{\nu}=D_eG_{\nu}(r,\phi)(s,\chi)\,\,\mbox{for }\,\,(s,\chi)\in\mathbb{R}^k\times C_n\,\,\mbox{and}\,\,\nu\in\{1,\ldots,n\}.
$$
The map
$$
{\mathcal U}\times\mathbb{R}^k\times C_n\ni(r,\phi,s,\chi)\mapsto D_eG(r,\phi)(s,\chi)\in\mathbb{R}^n
$$
is continuous, as required in property (He).

3. On $\Delta$. The chain rule shows that $\Delta:{\mathcal U} \to\mathbb{R}^k$ given by Eq. (9) is continuously differentiable. For every $(r,\phi)\in {\mathcal U}\subset I^k\times C^1_n$ and for all $\psi\in C^1_n$  we have
$$
D_2\Delta(r,\phi)\psi=D_2\delta(r,Q(r,\phi))D_2Q(r,\phi)\psi=D_2\delta(r,Q(r,\phi))Q(r,\psi)
$$
as each map $Q(r,\cdot)$, $r\in J^k$, is linear. The last expression also defines linear extensions $D_{2,e}\Delta(r,\phi):C_n\to\mathbb{R}^k$. Using continuity of $Q:J^k\times C_n\to F$ we obtain that the map
$$
{\mathcal U}\times C_n\ni(r,\phi,\chi)\mapsto D_{2,e}\Delta(r,\phi)\chi\in\mathbb{R}^k
$$
is continuous, as required in property (He). $\Box$

\begin{proposition}
For the system (3,4) with $G$ and $\Delta$ given by (8-11),\\
 $\emptyset\neq M_{G,\Delta}\subset {\mathcal U}$.
\end{proposition}

{\bf Proof.} By \cite[Proposition 2.3]{W5} there exists $\phi\in X_f$ for $f$ defined by $f(\phi)=g(v)$ where $v\in V$ has the components $v_{\iota}=\phi_j(-d_{\kappa}(L\phi))$ with 
$j\in\{1,\ldots,n\}$ and $\kappa\in\{1,\ldots,k\}$ determined by $(\kappa-1)n+j=\iota\in\{1,\ldots,kn\}$. In particular, $L\phi\in W$. Define $r\in [-h,0]^k\subset I^k$ by $r_{\kappa}=-d_{\kappa}(L\phi)\in[-h,0]$, $\kappa=1,\ldots,k$. It follows that
$Q(r,\phi)=L\phi\in W$, and for all $j\in\{1,\ldots,n\}$ and all $\kappa\in\{1,\ldots,k\}$, $E\phi_j(r_{\kappa})=\phi_j(r_{\kappa})=\phi_j(-d_{\kappa}(L\phi))$, which gives $\widehat{(r,\phi)}=v\in V$. We conclude that $(r,\phi)\in {\mathcal U}$. Moreover, by $\phi\in X_f$, $\phi'(0)=f(\phi)=g(v)=g(\widehat{(r,\phi)})=G(r,\phi)$. Also, $\Delta(r,\phi)=\delta(r,Q(r,\phi))=d(L\phi)+r=0$, with $d=(d_1,\ldots,d_k)$, and $\det\, D_1\Delta(r,\phi)=\det (id_{\mathbb{R}^k})=1\neq0$.
Altogether, $(r,\phi)\in M_{G,\Delta}.\quad\Box$

On vectorspaces of continuous linear maps from one normed vectorspace into another one, as well as on the vectorspace $\mathbb{R}^{n\times n}$ of $n\times n$-matrices with real entries, we use the norm given by $|A|=\sup_{|x|\le1}|Ax|$.

A product $C^1_n\times\mathbb{R}^n\to C^1_n$ which is used in the sequel and notated by $\phi\odot c$ is defined componentwise by $(\phi\odot c)_j=c_j\phi_j\in C^1,\quad j=1,\ldots,n.$

\section{Preparations for spaces complementary to $X_0$}

\begin{proposition}
Suppose $q:J^k\times C\to F$ is continuous with continuously differentiable restriction to $I^k\times C^1$, each map $q(r,\cdot):C\to F$, $r\in J^k$, is linear, and condition (Hq) is satisfied. Let $\epsilon>0$ be given. There exists a continuous map $\chi=\chi_q$ from $J^k$ into  $C^1$ with continuously differentiable restriction to $I^k$ such that for every $r\in J^k$ we have
$$
q(r,\chi(r))=0,\quad[\chi(r)]'(0)=1,\quad|\chi(r)|_C<\epsilon.
$$ 
\end{proposition}

{\bf Proof.} 1. For $q=0$ choose $\chi_q$ to be constant with its single value $\psi\in C^1$ satisfying $\psi'(0)=1$ and $|\psi|_C<\epsilon$.

2. The case $q\neq0$. 

2.1. Construction of a space $Y\subset C^1$ which is complementary in $C$ for each nullspace $q(r,\cdot)^{-1}(0)$, $r\in J^k$.  Recall the space $F_q$ and its dimension $d=d_q$, and the functions $\beta_1,\ldots,\beta_d$ in $C$ from property (Hq). For every $r\in J^k$ and all vectors $c$ in the unit sphere $S^{d-1}\subset\mathbb{R}^d$, $\sum_{m=1}^dc_mq(r,\beta_m)\neq0$. By continuity there are neighbourhoods $N_{r,c}$ of $(r,c)$ in $J^k\times S^{d-1}$ and $U_{r,c}$ of $(\beta_1,\ldots,\beta_d)$ in $C^d$ such that $\sum_{m=1}^d\tilde{c}_mq(\tilde{r},\tilde{\beta}_m)\neq0$ for all $(\tilde{r},\tilde{c})\in N_{r,c}$ and all $(\tilde{\beta}_1,\ldots,\tilde{\beta}_d)\in U_{r,c}$. By compactness $J^k\times S^{d-1}$ is covered by a finite collection of neighbourhoods $N_{r,c}$, say, by $N_{r_1,c_1},\ldots,N_{r_p,c_p}$ for $p\in\mathbb{N}$. It follows that for all $(r,c)\in J^k\times S^{d-1}$ and for all 
$$
(\tilde{\beta}_1,\ldots,\tilde{\beta}_d)\in\cap_{\pi=1}^p U_{r_{\pi},c_{\pi}}
$$
we have  $\sum_{m=1}^dc_mq(r,\tilde{\beta}_m)\neq0$. In $\cap_{\pi=1}^p U_{r_{\pi},c_{\pi}}$ we find $(\gamma_1,\dots,\gamma_d)\in (C^1)^d$ with $\gamma_1'(0)=0,\ldots,\gamma_d'(0)=0$. For all $r\in J^k$ and all $c\in S^{d-1}$,
 $\sum_{m=1}^dc_mq(r,\gamma_m)\neq0$, so $q(r,\gamma_1),\ldots,q(r,\gamma_d)$ are linearly independent and form a basis of $F_q$. According to a remark below condition (Hq) in Section 1 we obtain that for every $r\in J^k$,
$$
C=q(r,\cdot)^{-1}(0)\oplus\sum_{m=1}^d\mathbb{R}\,\gamma_m.
$$

2.2. Some bounds. Choose an isomorphism $\tau:F_q\to\mathbb{R}^d$. The maps
$$
J^k\ni r\mapsto\tau\,q(r,\gamma_m)\in\mathbb{R}^d,\quad m=1,\ldots,d,
$$
are continuous, as well as the matrix-valued map
$$
J^k\ni r\mapsto(\tau\,q(r,\gamma_1)\cdots\tau\,q(r,\gamma_d))^{-1}\in\mathbb{R}^{d\times d},
$$
whose restriction to the open set $I^k\subset\mathbb{R}^k$ is continuously differentiable. 
There is a constant $c_{\tau}\ge0$ with
$$
|(\tau\,q(r,\gamma_1)\cdots\tau\,q(r,\gamma_d))^{-1}|_{L_c(\mathbb{R}^d,\mathbb{R}^d)}\le c_{\tau}\quad\mbox{for all}\quad r\in J^k.
$$
As each set
$\{\tau\,q(r,\phi):r\in J^k\}$, $\phi\in C$, is bounded the Principle of Uniform Boundedness yields
$$
c_q=\sup_{r\in J^k}|\tau\,q(r,\cdot)|_{L_c(C,F)}<\infty.
$$  

2.3. Finding the desired map. Choose $\phi\in C^1$ with $\phi'(0)=1$ and
$$
|\phi|_C<\frac{\epsilon}{1+c_{\tau}c_q}
$$
and define $\chi(r)$ for $r\in J^k$  by projecting $\phi$ along $Y=\sum_{m=1}^d\mathbb{R}\gamma_m$ onto the nullspace
$q(r,\cdot)^{-1}(0)$, that is, $\chi(r)=\phi-\sum_{m=1}^dc_m(r)\gamma_m$ and 
$$
0=\tau\,q(r,\chi(r))=\tau\,q(r,\phi-\sum_{m=1}^dc_m(r)\gamma_m)=\tau\,q(r,\phi)-\sum_{m=1}^dc_m(r)\tau\,q(r,\gamma_m),
$$
or,
$$
c(r)=(c_1(r),\ldots,c_d(r))^{tr}=(\tau\,q(r,\gamma_1)\cdots\tau\,q(r,\gamma_d))^{-1}\cdot\tau\,q(r,\phi).
$$
Notice that $\chi$ is continuous, and that its restriction to $I^k$ is continuously differentiable, and that for every $r\in J$ we have
\begin{eqnarray}
[\chi(r)]'(0)=\phi'(0)-\sum_{m=1}^dc_m(r)\gamma_m'(0)=1,\nonumber\\
q(r,\chi(r))=0\quad\mbox{(by construction),\quad and}\nonumber\\
|\chi(r)|_C\le|\phi|_C+c_{\tau}c_q|\phi|_C<\epsilon.\quad\Box\nonumber
\end{eqnarray}

Proposition 2.1 may be viewed as a generalization of \cite[Proposition 2.1]{KW2}. However, the proof given above is different from the proof of \cite[Proposition 2.1]{KW2}. The next proposition corresponds to \cite[Proposition 2.2]{KW2}.
 
\begin{proposition}
Assume $q:J^k\times C\to F$ satisfies the hypothesis of Proposition 2.1. Let a continuous function $h_V:V\to(0,\infty)$ be given. There exists a continuously differentiable map $\chi=\chi_{q,V}$ from $I^k\times V$ into $C^1$ such that for all $r\in I^k$ and all $v\in V$ we have
$$
q(r,\chi(r,v))=0,\quad [\chi(r,v)] '(0)=1,\quad |\chi(r,v)|_C\le  h_V(v),
$$
and for every $\iota\in\{1,\ldots,kn\}$,
$$
|D_{k+\iota}\chi(r,v)1|_C\le h_V(v).
$$
\end{proposition}

Parts 1 and 2 of the proof of Proposition 2.2 can be taken from  the proof of \cite[Proposition 2.2]{KW2}. Nevertheless we prefer to include these parts here, in order to keep the presentation self-contained.

\medskip

{\bf Proof.} 1. There is a sequence of non-empty open subsets $V_{j1},V_{j2},V_j$ of $V$, with $j\in\mathbb{N}$, such that
$$
\bigcup_{j=1}^{\infty}V_j=V,
$$
and for every $j\in\mathbb{N}$,
$$
V_{j1}\subset\subset V_{j2}\subset\subset V_j\quad\mbox{and}\quad V_j\subset\subset V_{j+1,1}.
$$
With $\overline{V_{02}}=\emptyset$ we have that for each integer $j\ge 1$,
$$
(V_{j+1}\setminus\overline{V_{j2}})\cap(V_j\setminus\overline{V_{j-1,2}})=V_j\setminus\overline{V_{j2}}
$$
while for integers $j\ge1$ and $k\ge j+2$, 
$$
(V_k\setminus\overline{V_{k-1,2}})\cap(V_j\setminus\overline{V_{j-1,2}})=V_j\setminus\overline{V_{k-1,2}}\subset V_j\setminus \overline{V_{j+1,2}}=\emptyset. 
$$
 
2. For every $j\in\mathbb{N}$ choose a continuously differentiable function
$$
a_j:\mathbb{R}^{kn}\to[0,1]
$$
with
$$
a_j(v)=1\,\,\mbox{on}\,\,\overline{V_{j1}},\quad a_j(v)=0\,\,\mbox{on}\,\,\mathbb{R}^{kn}\setminus V_{j2}.
$$
For  every $j\in\mathbb{N}$ choose an upper bound
$$
A_j>1+\sum_{\mu=1}^{kn}\max_{v\in\mathbb{R}^{kn}}|D_{\mu}a_j(v)1|=\max_{v\in\overline{V_j}}|a_j(v)|+\sum_{\mu=1}^{kn}\max_{v\in\overline{V_j}}|D_{\mu}a_j(v)1|
$$
so that the sequence $(A_j)_1^{\infty}$ in $[1,\infty)$ is increasing.

\medskip

The sequence $(h_j)_1^{\infty}$ given by $h_j=\min_{v\in\overline{V_j}}h(v)>0$ is nonincreasing. We have 
$$
\frac{h_j}{2A_j}\le h_j\,\,\mbox{for all}\,\,j\in\mathbb{N},
$$
and the sequence $(h_j/2A_j)_{j=1}^{\infty}$ is decreasing.

3. A sequence of maps on the sets $I^k\times(V_j\setminus\overline{V_{j-1,2}})$. For $j\in\mathbb{N}$ set $\epsilon_j=\frac{h_j}{2A_j}$ and apply Proposition 2.1 to $q$ and $\epsilon=\epsilon_j$. This yields a sequence of continuously differentiable maps $\chi_j:I^k\to C^1$ such that for every $j\in\mathbb{N}$ and all $r\in I^k$,
$$
q(r,\chi_j(r))=0,\quad[\chi_j(r)]'(0)=1,\quad|\chi_j(r)|_C\le\epsilon_j.
$$
The maps
$$
H_j:I^k\times(V_j\setminus\overline{V_{j-1,2}})\to C^1,\quad j\in\mathbb{N},
$$
given by
$$
H_j(r,v)=a_j(v)\chi_j(r)+(1-a_j(v))\chi_{j+1}(r)
$$
are continuously differentiable. For all $j\in\mathbb{N}, r\in I^k,$ and $v\in V_j\setminus\overline{V_{j-1,2}}$ we have
\begin{eqnarray}
q(r,H_j(r,v)) & = & q(r,a_j(v)\chi_j(r)+(1-a_j(v))\chi_{j+1}(r))\nonumber\\
& = & a_j(v)q(r,\chi_j(r))+(1-a_j(v))q(r,\chi_{j+1}(r))=0,\nonumber\\
(H_j(r,v))'(0) & = & a_j(v)[\chi_j(r)]'(0)+(1-a_j(v))[\chi_{j+1}(r)]'(0)=1,\nonumber
\end{eqnarray}
and
\begin{eqnarray}
|H_j(r,v)|_C & \le &  a_j(v)|\chi_j(r)|_C+(1-a_j(v))|\chi_{j+1}(r)|_C\nonumber\\
& \le & a_j(v)\frac{h_j}{2A_j}+(1-a_j(v))\frac{h_{j+1}}{2A_{j+1}}\nonumber\\
& \le & a_j(v)\frac{h_j}{2A_j}+(1-a_j(v))\frac{h_j}{2A_j}\nonumber\\
& \le & h_j\quad\mbox{(with}\quad A_j\ge1)\nonumber\\
& \le &  h_V(v)\quad\mbox{(since}\quad v\in V_j).\nonumber
\end{eqnarray}
Moreover, for the same $j,r,v$ and for all $\iota\in\{1,\ldots,kn\}$,
$$
D_{k+\iota}H_j(r,v)1=D_{\iota}a_j(v)1\cdot\chi_j(r)-D_{\iota}a_j(v)1\cdot\chi_{j+1}(r),
$$
hence
\begin{eqnarray}
|D_{k+\iota}H_j(r,v)1|_C & \le & A_j\left(\frac{h_j}{2A_j}+\frac{h_{j+1}}{2A_{j+1}}\right)\nonumber\\
& \le & A_j\left(\frac{h_j}{2A_j}+\frac{h_j}{2A_j}\right)=h_j\le h_V(v)\quad\mbox{(since}\quad v\in V_j).\nonumber
\end{eqnarray}

4. It remains to show that the maps $H_j$ define a map on $I^k\times V$. This follows from Part 1 of the proof provided that for each $j\in\mathbb{N}$ the maps $H_j$ and $H_{j+1}$ coincide on the intersection  of their domains, which is
$$
(I^k\times(V_{j+1}\setminus\overline{V_{j,2}}))\cap (I^k\times(V_j\setminus\overline{V_{j-1,2}}))=I^k\times(V_j\setminus\overline{V_{j2}}).
$$
Proof of this: For $j\in\mathbb{N},r\in I^k,v\in V_j\setminus\overline{V_{j2}}$,
$$
H_{j+1}(r,v)=a_{j+1}(v)\chi_{j+1}(r)+(1-a_{j+1}(v))\chi_{j+2}(r)=\chi_{j+1}(r)
$$
since $a_{j+1}(v)=1$  on $\overline{V_{j+1,1}}\supset V_j\supset V_j\setminus\overline{V_{j2}}$, and
$$
H_j(r,v)=a_j(v)\chi_j(r)+(1-a_j(v))\chi_{j+1}(r)=\chi_{j+1}(r)
$$
since $a_j(v)=0$ on $\mathbb{R}^{kn}\setminus\overline{V_{j2}}\supset V_j\setminus\overline{V_{j2}}$. $\Box$

\section{Mapping the solution manifold into $\mathbb{R}^k\times X_0$}

The construction of the desired map begins with the choice of
$h_V=h_g$ where
$$
h_g(v)=\frac{\min\{1,dist(v,\mathbb{R}^{kn}\setminus V)\}}{6kn(1+\max_{\iota,j}|\partial_{\iota}g_j(v)|+\max_j|g_j(v)|)}
$$
in case $V\neq\mathbb{R}^{kn}$ and
$$
h_g(v)=\frac{1}{6kn(1+\max_{\iota,j}|\partial_{\iota}g_j(v)|+\max_j|g_j(v)|)}
$$
in case $V=\mathbb{R}^{kn}$. For $j\in\{1,\dots,n\}$ consider $q=Q_j$ where
$Q_j:J^k\times C\to F$ is given by $Q_j(r,\phi)=Q(r,\phi\cdot e_j)$.
Proposition 2.2 yields $n$ continuously differentiable maps
$\chi_{g,j}:I^k\times V\to C^1$, $j=1,\ldots,n$, which satisfy
$$
Q_j(r,\chi_{g,j}(r,v))=0,\quad[\chi_{g,j}(r,v)]'(0)=1,\quad\mbox{and}\quad|\chi_{g,j}(r,v)|_C\le h_g(v)
$$
for all $(r,v)\in I^k\times V$, and
$$
|\partial_{k+\nu}\chi_{g,j}(r,v)|_C\le h_g(v)\quad\mbox{for all}\quad\nu\in\{1,\ldots,kn\}\quad\mbox{and all}\quad(r,v)\in I^k\times V.
$$
For the 
continuously differentiable map $\chi_g:I^k\times V\to C^1_n$ defined by
$$
\chi_g(r,v)=(\chi_{g,1}(r,v),\ldots,\chi_{g,n}(r,v))
$$ 
we obtain that for all $c\in\mathbb{R}^n$ and for all $(r,v)\in I^k\times V$,
$$
\left[\sum_{j=1}^nc_j\chi_{g,j}(r,v)\cdot e_j\right]'(0)=\sum_{j=1}^nc_j[\chi_{g,j}(r,v)]'(0)e_j=\sum_{j=1}^nc_je_j=c.
$$
This implies that for every $(r,v)\in I^k\times V$ the elements
$\chi_ {g,j}(r,v)\cdot e_j$, $j=1,\ldots,n$, are linearly independent and that all spaces $\sum_{j=1}^n\mathbb{R}\chi_{g,j}(r,v)\cdot e_j$, $(r,v)\in I^k\times V$, are complementary to $X_0$ in $C^1_n$.
Notice also that 
$$
\sum_{j=1}^nc_j\chi_{g,j}(r,v)\cdot e_j=c\odot\chi_g(r,v)
$$
for all $c\in\mathbb{R}^n$ and all $(r,v)\in I^k\times V$.

Now consider the continuously differentiable maps
$A:{\mathcal U}\to C^1_n$ given by
$$
A(r,\phi)=\phi-g(v)\odot\chi_g(r,v),\quad v=\widehat{(r,\phi)},
$$
and
$T:{\mathcal U}\to\mathbb{R}^k\times C^1_n$ given by
$$
T(r,\phi)=(r,A(r,\phi)).
$$
Part (iii) of the next result shows in particular that $T$ maps the solution manifold 
$$
M_{G,\Delta}=\{(r,\phi)\in{\mathcal U}:\phi'(0)=G(r,\phi),\Delta(r,\phi)=0,D_1\Delta(r,\phi)\neq0\}
$$ 
into the subset of the strip $I^k\times X_0$ given by the 
`algebraic' equation $\Delta(r,\psi)=0$.
  
\begin{proposition}
(i) For all $c\in\mathbb{R}^n$, $r\in I^k$, and $v\in V$, 
$$
Q(r,c\odot\chi_g(r,v))=0.
$$

(ii) For all $r\in I^k$ and $v\in V$, $[\chi_g(r,v)]'(0)=(1,\dots,1)^{tr}$.

(iii) For every $(r,\phi)\in {\mathcal U}$,
$$
[A(r,\phi)]'(0)=\phi'(0)-g(v)\quad\mbox{with}\quad v=\widehat{(r,\phi)},
$$
and
$$
Q(r,A(r,\phi))=Q(r,\phi),\quad\Delta(T(r,\phi))=\Delta(r,\phi).
$$
In case $(r,\phi)\in M_{G,\Delta}$, 
$$
[A(r,\phi)]'(0)=0.
$$

(iv) For every $(r,\phi)\in {\mathcal U}$ with $(A(r,\phi))'(0)=0$ we have $$
\phi'(0)=g(\widehat{(r,\phi)})=G(r,\phi).
$$
If in addition $\Delta(r,\phi)=0$ and $\det\,D_1\Delta(r,\phi)\neq0$ then $(r,\phi)\in M_{G,\Delta}$.
\end{proposition}

{\bf Proof.} 1. On (i). 
\begin{eqnarray}
Q(r,c\odot\chi_g(r,v)) & = & Q(r,\sum_{j=1}^nc_j\chi_{g,j}(r,v)\cdot e_j)=\sum_{j=1}^nc_jQ(r,\chi_{g,j}(r,v)\cdot e_j)\nonumber\\
& = & \sum_{j=1}^nc_jQ_j(r,\chi_{g,j}(r,v))=\sum_{j=1}^nc_j0=0,\nonumber
\end{eqnarray}
due to Proposition 2.2.

2. On (ii). The components of $[\chi_g(r,v)]'(0)$ are $[\chi_{g,j}(r,v)]'(0)=1$, for $j=1,\ldots,n$.

3. On (iii). For every $(r,\phi)\in {\mathcal U}$, with $v=\widehat{(r,\phi)}$, 
\begin{eqnarray}
[A(r,\phi)]'(0) & = & \phi'(0)-\left[\sum_{j=1}^ng_j(v) \chi_{g,j}(r,v)\cdot e_j\right]'(0)\nonumber\\
& = & \phi'(0)-\sum_{j=1}^ng_j(v)[\chi_{g,j}(r,v)]'(0)e_j\nonumber\\ 
& = & \phi'(0)-\sum_{j=1}^ng_j(v)e_j=\phi'(0)-g(v),\nonumber
\end{eqnarray}
which implies $[A(r,\phi)]'(0)=0$ in case $(r,\phi)\in M_{G,\Delta}$.
From assertion (i),
$$
Q(r,A(r,\phi))=Q(r,\phi)-Q(r,g(v)\odot\chi_g(r,v))=Q(r,\phi)
$$
for every $(r,\phi)\in {\mathcal U}$. It follows that $Q(r,A(r,\phi))=Q(r,\phi)\in W$ and 
$$
\Delta(T(r,\phi))=\Delta(r,A(r,\phi))=\delta(r,Q(A(r,\phi)))=\delta(r,Q(r,\phi))=\Delta(r,\phi).
$$

4. On (iv). For $(r,\phi)\in {\mathcal U}$ with $(A(r,\phi))'(0)=0$ we obtain from Part (iii)
$0=\phi'(0)-g(v)$ with $v=\widehat{(r,\phi)}$. 
The remaining part of the assertion is now obvious from the definition of $M_{G,\Delta}$. $\Box$

\section{Representation by an algebraic equation}

In this section we show that the map $T$ is a diffeomorphism onto an open subset of $\mathbb{R}^k\times C^1_n$ and identify the image of the solution manifold. 
The proof that $T$ is invertible relies on the relation between $v=\widehat{(r,\phi)}$ for $(r,\phi)\in{\mathcal U}$ and $y=\widehat{(r,\psi)}$ for $(r,\psi)=T(r,\phi)=(r,A(r,\phi))$. For every $\iota\in\{1,\ldots,kn\}$, with $\kappa\in\{1,\ldots,k\}$ and $j\in\{1,\ldots,n\}$ given by $\iota=(\kappa-1)n+j$, we have
\begin{eqnarray}
y_{\iota} & = & (E\psi_j)(r_{\kappa})=(E(A(r,\phi)_j))(r_{\kappa})\nonumber\\
& = & (E\phi_j)(r_{\kappa})-(E(g(v)\odot\chi(r,v))_j)(r_{\kappa})\nonumber\\
& = & v_{\iota}-(E(g_j(v)\chi_{g,j}(r,v)))(r_{\kappa})\nonumber\\
& = & v_{\iota}-g_j(v)(E(\chi_{g,j}(r,v)))(r_{\kappa})\nonumber\\
& = & v_{\iota}-g_j(v)ev(E(\chi_{g,j}(r,v)),r_{\kappa}),\nonumber
\end{eqnarray}
hence
$$
y=S(r,v)=v-R(r,v)
$$
with the continuously differentiable 
$$
S:I^k\times V\to\mathbb{R}^{kn}\quad\mbox{and}\quad
R:I^k\times V\to\mathbb{R}^{kn}
$$
given by  
$$
R_{\iota}(r,v)=g_j(v)ev(E(\chi_{g,j}(r,v)),r_{\kappa})\quad\mbox{for}\quad\iota=1,\ldots,kn,
$$ 
with  $\kappa\in\{1,\ldots,k\}$ and $j\in\{1,\ldots,n\}$ determined by $\iota=(\kappa-1)n+j$.

The next result uses the properties of the maps $\chi_{g,j}$ in order to show that the map $S$ is a perturbation of the identity which is under control by a Lipschitz estimate and smallness close to the boundary of its domain.

\begin{proposition}
(i) For every $(r,v)\in I^k\times V$, 
$$
|D_2R(r,v)|_{L_c(\mathbb{R}^{kn},\mathbb{R}^{kn})}\le\frac{1}{2}.
$$
(ii) In case $V\neq\mathbb{R}^{kn}$,
$$
|R(r,v)|\le\frac{1}{2}dist(v,\mathbb{R}^{kn}\setminus V)\quad\mbox{for every}\quad(r,v)\in I^k\times V.
$$
\end{proposition}

{\bf Proof.} 1. On (i). For $r\in I^k$ define $R_r:V\to\mathbb{R}^{kn}$ by $R_r(v)=R(r,v)$. Then, for every $v\in V$ and all $w\in\mathbb{R}^{kn}$ with $|w|\le1$,
\begin{eqnarray}
|D_2R(r,v)w| & = & |DR_r(v)w|=\sqrt{\sum_{\iota=1}^{kn}\left(\sum_{\nu=1}^{kn}\partial_{\nu}R_{r,\iota}(v)w_{\nu}\right)^2}\nonumber\\
& \le &   \sqrt{\sum_{\iota=1}^{kn}\left(\sum_{\nu=1}^{kn}(\partial_{\nu}R_{r,\iota}(v))^2\right)\left(\sum_{\nu=1}^{kn}w_{\nu}^2\right)}\nonumber\\    & & \mbox{(by the Cauchy-Schwartz inequality)}\nonumber\\
& \le & \sqrt{kn\cdot kn\,\max_{\nu,\iota}(\partial_{\nu}R_{r,\iota}(v))^2\cdot 1}\nonumber\\
& = & kn\,\max_{\nu,\iota}|\partial_{\nu}R_{r,\iota}(v)|,\nonumber
\end{eqnarray}
hence
$$
|D_2R(r,v)|_{L_c(\mathbb{R}^{kn},\mathbb{R}^{kn})}\le  kn\,\max_{\nu,\iota}|\partial_{\nu}R_{r,\iota}(v)|.
$$
For each $v\in V$, $\iota\in\{1,\ldots,kn\}$, and $\kappa\in\{1,\ldots,k\}$ and $j\in\{1,\ldots,n\}$ defined by $\iota=(\kappa-1)n+j$, we have
$$
R_{r,\iota}(v)=g_j(v)ev(E(\chi_{g,j}(r,v)),r_{\kappa}).
$$
Partial differentiation for every $\nu\in\{1,\ldots,kn\}$ yields
$$
\partial_{\nu}R_{r,\iota}(v)=ev(E(\chi_{g,j}(r,v)),r_{\kappa})\partial_{\nu}g_j(v)+g_j(v)[E(\partial_{k+\nu}\chi_{g,j}(r,v))](r_{\kappa}),
$$
by means of the product rule, by the formula $Dev_1(\phi,s)(\tilde{\phi},\tilde{s})=ev(\tilde{\phi},s)+\tilde{s}\phi'(s)$, by linearity of the continuous map $C^1\stackrel{E}{\to} C^1(J,\mathbb{R})$, and by the chain rule. It follows that
\begin{eqnarray}
|\partial_{\nu}R_{r,\iota}(v)| & \le & |\partial_{\nu}g_j(v)||E(\chi_{g,j}(r,v))|_C+|g_j(v)||E(\partial_{k+\nu}\chi_{g,j}(r,v))|_C\nonumber\\
& \le & 3(|\partial_{\nu}g_j(v)||\chi_{g,j}(r,v)|_C+|g_j(v)||\partial_{k+\nu}\chi_{g,j}(r,v)|_C)\nonumber\\
& \le & 3(|\partial_{\nu}g_j(v)|+|g_j(v)|)h_g(v),\nonumber
\end{eqnarray}
by the construction of $\chi_{g,j}$. Hence
\begin{eqnarray}
|D_2R(r,v)|_{L_c(\mathbb{R}^{kn},\mathbb{R}^{kn})} & \le &  kn\,\max_{\nu,\iota}|\partial_{\nu}R_{r,\iota}(v)|\nonumber\\
& \le & 3kn(\max_{\nu,j})|\partial_{\nu}g_j(v)|+\max_j|g_j(v)|)h_g(v)\le\frac{1}{2}.\nonumber
\end{eqnarray}

2. On (ii). In case $V\neq\mathbb{R}^{kn}$, for $r\in I^k$, $v\in V$, $\iota\in\{1,\ldots,kn\}$, and $\kappa\in\{1,\ldots,k\}$ and $j\in\{1,\ldots,n\}$ defined by $\iota=(\kappa-1)n+j$, we have
\begin{eqnarray}
|R_{\iota}(r,v)| & = &|g_j(v)||(E(\chi_{g,j}(r,v))(r_{\kappa})|\nonumber\\
& \le & \max_j|g_j(v)||E(\chi_{g,j}(r,v))|_C\le 3\,\max_j|g_j(v)||\chi_{g,j}(r,v)|_C\nonumber\\
& \le & 3\,\max_j|g_j(v)||h_g(v)|,\nonumber
\end{eqnarray}
by the construction of $\chi_{g,j}$. Finally,
$$
|R(r,v)|\le\sqrt{kn}\max_{\iota}|R_{\iota}(r,v)|\le3\sqrt{kn}\max_j|g_j(v)||h_g(v)|\le\frac{1}{2}\,dist(v,\mathbb{R}^{kn}\setminus V).
$$
$\Box$

For $r\in I^k$ define $S_r:V\to\mathbb{R}^{kn}$ by $S_r(v)=S(r,v)$.
As a consequence of Proposition 4.1 we have

\begin{proposition}
(i) The set 
$$\cup_{r\in I^k}\{r\}\times S_r(V)\subset I^k\times\mathbb{R}^{kn}
$$ 
is open.

(ii) Each map $S_r$, $r\in I^k$, is a diffeomorphism onto the open set $$
S_r(V)=S(\{r\}\times V)\subset\mathbb{R}^{kn}.
$$

(iii) The map
$$
\cup_{r\in I^k}\{r\}\times S_r(V)\ni(r,y)\mapsto S_r^{-1}(y)\in V\subset\mathbb{R}^{kn}
$$
is continuously differentiable.
\end{proposition}

For the proof, see \cite[Proposition 3.2]{KW2}. 
We proceed to the verification that $T$ is invertible. As the sets $W$ and $\cup_{r\in I^k}\{r\}\times S_r(V)$ are open and the maps $I^k\times C_n\ni(r,\psi)\mapsto \widehat{(r,\psi)})\in\mathbb{R}^{kn}$ and $Q$ are continuous the set
\begin{eqnarray}
{\mathcal O} & = & \{(r,\psi)\in I^k\times\mathbb{R}^{kn}:Q(r,\psi)\in W,(r,\widehat{(r,\psi)})\in\cup_{r\in I^k}\{r\}\times S_r(V)\}\nonumber\\
& = & \{(r,\psi)\in I^k\times\mathbb{R}^{kn}:Q(r,\psi)\in W,\widehat{(r,\psi)})\in S_r(V)\}\nonumber
\end{eqnarray}
is open. The maps
$$
B:{\mathcal O}\to C^1_n,\quad B(r,\psi)=\psi+g(v)\odot\chi_g(r,v)\quad\mbox{with}\quad v=S_r^{-1}(y),y=\widehat{(r,\psi)},
$$
and
$$
Y:{\mathcal O}\to I^k\times C^1_n,\quad Y(r,\psi)=(r,B(r,\psi)),
$$
are continuously differentiable.

\begin{proposition}
(i) For every $(r,\psi)\in{\mathcal O}$, $Q(Y(r,\psi)=Q(r,\psi)$.

(ii) $T({\mathcal U})\subset{\mathcal O}$, and for every $(r,\phi)\in {\mathcal U}$, $Y(T(r,\phi))=(r,\phi)$.

(iii) $Y({\mathcal O})\subset {\mathcal U}$, and for every $(r,\psi)\in {\mathcal O}$, $T(Y(r,\psi))=(r,\psi)$.
\end{proposition}

{\bf Proof.} 1. On (i). For $(r,\psi)\in{\mathcal O}$ set $v=S_r^{-1}(y)$ with $y=\widehat{(r,\psi)}$. Using Proposition 3.1 (i) we infer
$$
Q(Y(r,\psi))=Q(r,B(r,\psi))=Q(r,\psi)+Q(r,g(v)\odot\chi_g(r,v))=Q(r,\psi).
$$

2. On (ii). Let $(r,\phi)\in {\mathcal U}$ be given and set $(r,\psi)=T(r,\phi)=(r,A(r,\phi))$.
By Proposition 3.1 (iii), $Q(r,\psi)=Q(r,A(r,\phi))=Q(r,\phi)\in W$.
Let $v=\widehat{(r,\phi)}$ and $y=\widehat{(r,\psi)}$. By the remarks at the begin of this section,  $\widehat{(r,\psi)}=y=S(r,v)=S_r(v)\in S_r(V)$. It follows that $(r,\psi)\in{\mathcal O}$, and we obtain $T({\mathcal U})\subset{\mathcal O}$. Moreover,  $Y(T(r,\phi))=Y(r,\psi)=(r,B(r,\psi))$ with
$$
B(r,\psi)=\psi+g(v)\odot\chi_g(r,v)=A(r,\phi)+g(v)\odot\chi_g(r,v)=\phi,$$
hence $Y(T(r,\phi))=(r,\phi)$.\\

3. On (iii). 

3.1. Proof of $Y({\mathcal O})\subset {\mathcal U}$. Let $(r,\psi)\in {\mathcal O}$ be given. Set $(r,\phi)=Y(r,\psi)=(r,B(r,\psi))$. By (i), $Q(r,\phi)=Q(r,\psi)\in W$. For $(r,\phi)\in {\mathcal U}$ it remains to show $\widehat{(r,\phi)}\in V$. In order to see this, recall that due to $(r,\psi)\in{\mathcal O}$ there exists $v\in V$ with  $y=\widehat{(r,\psi)}=S_r(v)$. We complete the proof by computing
$\widehat{(r,\phi)}=v$ : For $\iota\in\{1,\ldots,kn\}$ and $\kappa\in\{1,\ldots,k\}$, $j\in\{1,\ldots,n\}$ given by $\iota=(\kappa-1)n+j$,
\begin{eqnarray}
\widehat{(r,\phi)}_{\iota} & = & (E\phi_j)(r_{\kappa})= (E(B(r,\psi)_j))(r_{\kappa})\nonumber\\
& = & (E\psi_j)(r_{\kappa})+[E(g_j(v)\chi_{g,j}(r,v))](r_{\kappa})\nonumber\\
& = & \widehat{(r,\psi)}_{\iota}+g_j(v)[E(\chi_{g,j}(r,v))](r_{\kappa})=y_{\iota}+R_{\iota}(r,v)\nonumber\\
& = & S_{\iota}(r,v)+R_{\iota}(r,v)=v_{\iota}.\nonumber
\end{eqnarray}

3.2. Proof of $T(Y(r,\psi))=(r,\psi)$ for $(r,\psi)\in{\mathcal O}$. Set $(r,\phi)=Y(r,\psi)$. Then
$$
\phi=B(r,\psi)=\psi+g(v)\odot \chi_g(r,v)\quad\mbox{with}\quad v=S_r^{-1}(y),\quad y=\widehat{(r,\psi}).
$$
In Part 3.1 we saw that $v=\widehat{(r,\phi)}$. Using this and $Y(r,\psi)\in {\mathcal U}$ we get $T(Y(r,\psi))=T(r,\phi)=(r,A(r,\phi))$ with
$$
A(r,\phi)=\phi-g(v)\odot\chi_g(r,v)=B(r,\psi)-g(v)\odot\chi_g(r,v)=\psi.
$$
Hence $T(Y(r,\psi))=(r,\psi)$. $\Box$

According to Proposition 4.3 the map $T$ is a diffeomorphism onto the open set ${\mathcal O}\subset I^k\times C^1_n$, and $T^{-1}=Y$. Recall the trivial solution manifold $X_0\subset C^1_n$ given by the equation $\chi'(0)=0$.

\begin{corollary}
\begin{eqnarray}
T(M_{G,\Delta}) & = & \{(r,\psi)\in{\mathcal O}:\psi'(0)=0,\Delta(r,\psi)=0,\det\,D_1\Delta(T^{-1}(r,\psi))\neq0\}\nonumber\\
& = & \{(r,\psi)\in(I^k\times X_0)\cap {\mathcal O}:\Delta(r,\psi)=0,\det\,D_1\Delta(T^{-1}(r,\psi))\neq0\}.\nonumber
\end{eqnarray}
\end{corollary}

{\bf Proof} of the first equation of the corollary: The definition of $M_{G,\Delta}$ in combination with $T({\mathcal U})\subset{\mathcal O}$ and Proposition 3.1 (iii) yields the inclusion
$$
T(M_{G,\Delta})\subset \{(r,\psi)\in{\mathcal O}:\psi'(0)=0,\Delta(r,\psi)=0,\det\,D_1\Delta(T^{-1}(r,\psi))\neq0\}.
$$

Conversely, let $(r,\psi)\in{\mathcal O}$ be given with
$\psi'(0)=0$, $\Delta(r,\psi)=0$, and $\det\,D_1\Delta(T^{-1}(r,\psi))\neq0$. Set $(r,\phi)=Y(r,\psi)=(r,B(r,\psi))$. Then $(r,\phi)\in {\mathcal U}$ (see Proposition 4.3 (iii)) and $(r,\psi)=T(r,\phi)$. Using 
$(r,A(r,\phi))=T(r,\phi)=(r,\psi)$ and $\psi'(0)=0$ we infer $[A(r,\phi)]'(0)=0$. Proposition 3.1 (iv) yields 
$\phi'(0)=g(\widehat{(r,\phi)})$. By Proposition 3.1 (iii),
$$
0=\Delta(r,\psi)=\Delta(T(r,\phi))=\Delta(r,\phi).
$$
Moreover, $\det D_1\Delta(r,\phi)=\det D_1\Delta(T^{-1}(r,\psi))\neq0.$
Altogether, $(r,\phi)\in M_{G,\Delta}$, and $(r,\psi)=T(r,\phi)\in T(M_{G,\Delta})$. $\Box$

\section{A graph representation for explicit  delays}

Recall Eq. (2) from Section 1, with maps $g:\mathbb{R}^{kn}\supset V\to\mathbb{R}^{kn}$,
$L:C_n\to F$, and $d_{\kappa}:F\supset W\to[0,h]$, $\kappa\in\{1,\ldots,k\}$. Let $d:W\to\mathbb{R}^k$ be the map with components $d_{\kappa}$, $\kappa=1,\ldots,k$. Eq. (2) is equivalent to the system (3,4) with (8) $G(r,\phi)=g(\widehat{(r,\phi)})$, (9) $\Delta(r,\phi)=\delta(r,Q(r,\phi))$,
\begin{equation}
Q(r,\phi)=L\phi\quad\mbox{for all}\quad(r,\phi)\in J^k\times C_n, 
\end{equation}
and
\begin{equation}
\delta(r,w)=d(w)+r\quad\mbox{for all}\quad(r,w)\in J^k\times W, 
\end{equation}
in the sense that the continuously differentiable solutions $[-h,t_e)\to\mathbb{R}^n$ of Eq. (2) coincide with the first components $[-h,t_e)\to\mathbb{R}^n$ of solutions to (3,4) specified by (8-11). 

\begin{corollary}
For the system (3,4) with $G$ and $\Delta$ given by (8-11),
$$
T(M_{G,\Delta})=\{(r,\psi)\in(I^k\times X_0)\cap {\mathcal O}:r=-d(L\psi)\}
$$	
\end{corollary}

{\bf Proof.} Observe that for the systems considered we have
$\det\,D_1\Delta(r,\phi)=1\neq0$ for all
$(r,\phi)\in I^k\times C_n$. Use the second equation in Corollary 4.4 in combination with the choice of $Q$ and $\delta$. $\Box$
 
\section{Appendix}

There are only 4 minor changes to be made in \cite{W4} in order to obtain a result for the solution manifold $M=M_{G,\Delta}$ 
of the system (3,4) with $G$ and $\Delta$ defined on a subset ${\mathcal U}$ of the open strip $ I^k\times C^1_n$  (instead of $U\subset(-h,0)^k\times C^1_n$ as originally in \cite{W4}), under the smoothness hypotheses for $G$ and $\Delta$  stated in Section 1. These changes are the following ones:
\begin{itemize}
\item At the begin of \cite[Section 2]{W4} replace {\it neighbourhood $V_0$ in $(-h,0)^k$} by
$$
 neighbourhood\,\,V_0\,\, in\,\,\mathbb{R}^k
$$
\item In the proof of \cite[Proposition 3.4]{W4} replace the definition of $S_x$ by
$$
S_x=\{(s,t)\in\mathbb{R}^k\times(0,t_e):(s,x_t)\in U\}.
$$
\item In the proof of \cite[Proposition 3.5]{W4} replace {\it ... open neighbourhoods $N_0\subset C^1$ of $\phi_0$ and $V_0\subset(-h,0)^ k$ of $s_0$ with ...} by
$$
open\,\, neighbourhoods\,\,N_0\subset C^1\,\,of\,\,\phi_0\,\,and\,\,V_0\subset\mathbb{R}^k\,\,of\,\,s_0\,\,with\,\, ...
$$
\item In the proof of \cite[Proposition 3.5]{W4} replace {\it ... consider $\rho:[0,t_e)\ni t\mapsto r(t_x+t)\in(-h,0)^k$...} by
$$
 ...\,\, consider\,\,\rho:[0,t_e)\ni t\mapsto r(t_x+t)\in\mathbb{R}^k\,\, ...
$$
\end{itemize}


\end{document}